\documentclass[11pt]{article}

\usepackage{amssymb,amsbsy,amsthm,amsmath,amsfonts}

\def\sh{\hbox {\rm sinh\,}}

\newcommand{\rr}{\stackrel {d}{=}}

\renewcommand{\Re}{{\rm I\kern-0.16em R}}

\def\@begintheorem#1#2{\trivlist \item[\hskip \labelsep{\bf #1\ #2}]}
\def\@opargbegintheorem#1#2#3{\trivlist
      \item[\hskip \labelsep{\bf #1\ #2\ (#3)}]}

\newtheorem{proposition}{Proposition}%[section] ger numering sektionsvis %denna med i latex

\newtheorem{thm}[proposition]{Theorem}

\newtheorem{corollary}[proposition]{Corollary}

\newtheorem{remark}[proposition]{Remark}

\def\sh{\hbox {\rm sinh}}
\def\P{{\bf P}}
\def\R{{\bf R}}
\def\R{{\bf R}}
\def\Q{{\bf Q}}
\def\E{{\bf E}}
\def\ep{{\epsilon}}

\def\cF{{\cal F}}
\def\cD{{\cal D}}

\def\cC{{\cal C}}
\def\cG{{\cal G}}

\def\ep{\varepsilon}

\begin{document}

\author{
Paavo Salminen\\{\small \AA bo Akademi University,}
	\\{\small Mathematical Department,}
	\\{\small F\"anriksgatan 3 B,}
	\\{\small FIN-20500 \AA bo, Finland,} 
	\\{\small \tt phsalmin@abo.fi}
\and
Pierre Vallois
\\{\small Universit\'e Henri Poincar\'e}
\\{\small D\'epartement de Math\'ematique}
\\{\small F-54506 Vandoeuvre les Nancy, France}
\\{\small \tt vallois@iecn.u-nancy.fr}
%{\small email: }
\and
Marc Yor\\
	{\small Universit\'e Pierre et Marie Curie,}\\
	{\small Laboratoire de Probabilit\'es }\\
	{\small et Mod\`eles al\'eatoires,}\\
	{\small 4, Place Jussieu, Case 188}\\
	{\small F-75252 Paris Cedex 05, France}
}
\vskip5cm

\title{
	On the excursion theory for linear diffusions}
%\date{} % den har med om datumet inte skall med

\maketitle

\begin{abstract} We present a number of important identities related
  to the excursion theory of linear diffusions. In particular, excursions straddling
an independent exponential time are studied in detail. Letting the parameter of 
the exponential time tend to zero it is seen that these results connect 
to the corresponding results for excursions of stationary diffusions (in stationary state). We
characterize also the laws of the diffusion prior and posterior to the last zero 
before the exponential time. It is proved using Krein's
representations that, e.g., the law of the length of the excursion
straddling an exponential time is infinitely divisible. As an illustration of the results 
we discuss Ornstein-Uhlenbeck processes.  
	\\ \\%\bigskip\noindent
	\\ \\ %\bigskip\noindent
\end{abstract}

\section{Introduction and preliminaries}
\label{sec1}

{\bf 1.1}\hskip2mm Throughout this paper, we shall assume that $X$ is a linear regular recurrent 
diffusion taking values in $\R_+$ with 0 an instantaneously reflecting boundary. Let $\P_x$ and $\E_x$ 
denote, respectively, the probability measure and the
expectation associated with $X$ when started from $x\geq 0.$ We assume that $X$ is 
defined in the canonical space $C$ of continuous functions $\omega:\R_+\mapsto \R_+.$ 
Let 
$$
{\cal C}_t:=\sigma\{\omega(s): s\leq t\}
$$  
denote the smallest $\sigma$-algebra making the co-ordinate mappings up to time $t$ measurable 
and take ${\cal C}$ to be the smallest $\sigma$-algebra including all $\sigma$-algebras ${\cal C}_t,\ t\geq 0.$

The excursion space for excursions from 0 to 0 associated with $X$ 
is a subset of $C,$ denoted by $E,$ and given by
\begin{eqnarray*}
&&\hskip-1cm
E:=\{\varepsilon\in C: \varepsilon(0)=0,\ \exists\ \zeta(\ep)>0\ 
{\rm such\ that\ } \ep(t)>0\ \forall\ t\in(0,\zeta(\ep))
\\ 
&&\hskip4cm
{\rm and}\ \ep(t)=0\ \forall\ t\geq \zeta(\ep)\}.
\end{eqnarray*}   
The notation  ${\cal E}_t$ is used for the trace of ${\cal C}_t$ on $E.$

As indicated in the title of the paper our aim is to gather a number of fundamental results 
concerning 
the excursion theory for the diffusion $X.$ In Section 2 the classical
descriptions, the first one due to It\^o and McKean and the second one
due to Williams, are presented. In Section 3 the stationary excursions
are discussed and, in particular, the description due to Bismut is reviewed.   
After this, in Section 4, we proceed by analyzing 
excursions straddling an exponential time.  The paper is concluded with an example on 
Ornstein-Uhlenbeck processes.

Our motivation for this work arose from different origins:   
%\begin{description}
\begin{itemize}
\item First, we would like to contribute to Professor It\^o's being awarded the 1st Gauss prize, 
by offering some discussion and illustration of K It\^o's excursion theory, see \cite{ito70}, when 
specialized to linear diffusions. The present paper also illustrates Pitman and Yor's discussion
 (see \cite{pitmanyor06} in this volume) of K. It\^o's general theory of excursions for a Markov process.
\item In the literature there seems to be lacking a
  detailed 
discussion on the excursion theory of linear diffusions. 
Information 
available has a very scattered character, see, e.g., Williams \cite{williams74}, 
Walsh \cite{walsh78b},  Pitman and Yor \cite{pitmanyor82}, \cite{pitmanyor96}, \cite{pitmanyor97},  \cite{pitmanyor99}, 
Rogers \cite{rogers89}, Salminen \cite{salminen97}. The general theory of excursions has been developed in 
It\^o \cite{ito70}, Meyer \cite{meyer71}, Getoor \cite{getoor79a}, Getoor and Sharpe \cite{getoorsharpe73a}, 
\cite{getoorsharpe73b}, \cite{getoorsharpe79}, \cite{getoorsharpe82}, Blumenthal \cite{blumenthal92}. 
Although the case with Brownian motion is well 
studied and understood, for textbook treatments see,  e.g., Revuz and Yor \cite{revuzyor01} 
and Rogers and Williams \cite{rogerswilliams87}, 
we find it important to highlight the main formulas for more general diffusions using the traditional Fellerian terminology and language. 
%of diffusions.% (see It\^o and McKean \cite{itomckean74}).
\item To generalize some recent results (see Winkel \cite{winkel05} 
and Bertoin, Fujita, Roynette and Yor \cite{bertoinetal06}) on infinite divisibility of the distribution 
of the length of the excursion of a diffusion straddling 
an independent exponential time.
\item The Ornstein-Uhlenbeck process is one of the most essential diffusions. To present in detail 
formulae for its excursions is 
important per se. One of the key tools hereby is the distribution 
of the first hitting time $H_y$ of the point $y$ from which the excursions are 
observed.    
For $y=0$ this distribution can be derived via Doob's transform (see Doob \cite{doob42}) 
which connects the Ornstein-Uhlenbeck process  
with standard Brownian motion (see Sato \cite{sato77}, and G\"oing-Jaeschke and Yor  \cite{going-jaeschkeyor03}). 
For arbitary $y$ the distribution is 
very complicated; for explicit expressions via series expansions, see  
Ricciardi and Sato \cite{ricciardisato88}, Linetsky \cite{linetsky04} 
and  Alili, Patie and Pedersen \cite{alilipatiepedersen05}. 
We will focus on excursions from 0 to 0 and relate our work to earlier papers 
by Hawkes and Truman \cite{hawkestruman91}, Pitman and Yor 
\cite{pitmanyor97}, and  Salminen \cite{salminen84b}. Due to the
symmetry of the Ornstein-Uhlenbeck process around 0, it is sufficient
for our purposes to consider only positive excursions - the treatment of negative ones
is similar - and view the process with values in $\R_+$ and 0 being
a reflecting boundary. 
%XXXX a real ``comedy of errors''XXXX
\end{itemize}
%\end{description}

\noindent
{\bf 1.2}\hskip2mm In this subsection we introduce the basic notation and facts concerning linear diffusions needed 
in the sequel. A main source of information remains It\^o and McKean \cite{itomckean74}, see also Rogers and 
Williams \cite{rogerswilliams87}, and Borodin and 
Salminen \cite{borodinsalminen02}.
\begin{description}
\item{(i)}\hskip3mm Speed measure $m$ associated with $X$ is a  measure on $\R_+$ which satisfies
for all $0<a<b<\infty$
$$
0<m((a,b))<\infty.
$$
For simplicity, it is assumed that $m$ 
does not have atoms. An important fact is that $X$ has a jointly continuous transition density 
$p(t;x,y)$ with respect to  $m,$ i.e., 
$$
\P_x(X_t\in A)=\int_A p(t;x,y)\, m(dy),
$$
where $A$ is a Borel subset of $\R_+.$ Moreover, $p$ is symmetric in $x$ and $y,$ that is,
$p(t;x,y)=p(t;y,x).$ The Green or the resolvent kernel of $X$ is defined for $\lambda>0$ as
$$
 R_\lambda(x,y)=\int_0^\infty dt\ {\rm e}^{-\lambda t}\,p(t;x,y).
 $$
\item{(ii)} Scale function $S$ is an increasing and continuous function which can be defined via the identity
\begin{equation}
\label{scale}
\P_x(H_a<H_b)=\frac{S(b)-S(x)}{S(b)-S(a)},\quad 0\leq a<x<b, 
\end{equation}
where $H_.$ denotes the first hitting time, i.e., 
$$
H_y:=\inf\{t:\, X_t=y\},\quad y\geq 0. 
$$
We normalize by setting $S(0)=0.$ Due to the 
recurrence assumption it holds $S(+\infty)=+\infty.$ 
Recall that $\{S(X_{t\wedge H_0}): t\geq 0\}$ is a continuous 
local $\P_x$-martingale for every $x\geq 0$ (see, e.g., 
Rogers and Williams \cite{rogerswilliams87} p. 276). It is easily proved that $S(X)=\{S(X_{t}): t\geq 0\}$ is a (recurrent) diffusion 
taking values in $\R_+.$  The scale function associated with $S(X)$ is the 
identity mapping $x\mapsto x,\, x\geq 0,$ and we say that $S(X)$ is in natural scale. Clearly, also for $S(X)$ the boundary point  
0 is instantaneously reflecting. Using the Skorokhod reflection equation it is seen that $S(X)$ is a $\P_x$-submartingale (cf. Meleard 
\cite{meleard86} Proposition 1.4 where the semimartingale decomposition is given in case there are two reflecting boundaries).  
\item{(iii)} The infinitesimal generator of $X$ can be expressed  as the generalized differential operator 
$$
\cG = \frac{d}{dm}\,\frac{d}{dS}
$$
acting  on functions $f$ belonging to the appropriately defined domain $\cD(\cG)$ of $\cG$ (see  It\^o and McKean \cite{itomckean74}, 
Freedman \cite{freedman71}, Borodin and Salminen \cite{borodinsalminen02}). In particular, since 0 is assumed to be reflecting 
then $f\in\cD(\cG)$ implies that 
$$
f^+(0):=\lim_{x\uparrow 0} \frac{f(x)-f(0)}{S(x)-S(0)}=0.
$$
\item{(iv)} The distribution of the first hitting time of a point $y>0$ has a $\P_x$-density:
$$
\P_x(H_y\in dt)=f_{xy}(t)\,dt.
$$ 
This density can be connected with the derivative of a transition density of a killed diffusion obtained from $X.$  
To explain this, introduce the sample paths 
$$
\widehat X^{(y)}_t:=
\begin{cases}
X_t, & t<H_y,\\
\partial,& t\geq H_y,
\end{cases}
$$
where $\partial$  is a point isolated from $\R_+$ (a "cemetary" point). Then $\{\widehat X^{(y)}_t:t\geq 0\}$ is a diffusion    
with the same scale and speed as $X.$ Let $\hat p$ denote the transition density of $\widehat X^{(y)}$ with respect to $m.$ 
Then, e.g., for $x>y$
\begin{equation}
\label{f00}
f_{xy}(t)=\lim_{z\downarrow y}\frac{\hat p(t;x,z)}{S(z)-S(y)}.
\end{equation}
For a fixed $x$ and $y,$ the mapping $t\mapsto f_{xy}(t)$ is continuous, 
as follows from the eigen-differential expansions and discussion
in It\^o and\break McKean p. 153 and 217 (see also Kent \cite{kent80}, \cite{kent82}). Recall also the 
following formula for the Laplace transform of $H_y$
\begin{equation}
\label{green}
\E_x\left({\rm e}^{-\alpha H_y}\right)
=\frac {R_\alpha(x,y)}{R_\alpha(y,y)},
\end{equation}
which leads to 
$$
\int_0^\infty m(dx)\,
\E_x\left({\rm e}^{-\alpha H_y}\right)
=\frac 1{\alpha R_\alpha(y,y)}.
$$   
\item{(v)} There exists a jointly continuous family of local times 
$$
\{ L^{(y)}_t: t\geq 0, y\geq 0\}
$$ 
such that  $X$ satisfies the occupation time formula
$$
\int_0^t ds \,h(X_s) =\int_0^\infty h(y) L^{(y)}_t m(dy),
$$
where $h$ is a nonnegative measurable function (see, e.g., Rogers and Williams \cite{rogerswilliams87} 
49.1 Theorem p. 289). Consequently,
$$
L^{(y)}_t=\lim_{\delta\downarrow 0}\frac 1{m((y-\delta,y+\delta))} \int_0^t 
{\bf 1}_{(y-\delta,y+\delta)}(X_s)\, ds.
$$
For a fixed $y$ introduce the inverse of $L^{(y)}$ via
$$
\tau^{(y)}_\ell:=\inf\{s: L^{(y)}_s>\ell\}.
$$
Then $\tau^{(y)}=\{\tau^{(y)}_\ell: \ell\geq 0\}$ 
is an increasing L\'evy process, in other words, a subordinator
and its L\'evy exponent is given by
\begin{eqnarray}
\label{e1}
&&\hskip-1cm
\nonumber
\E_y\left(\exp(-\lambda \tau^{(y)}_\ell)\right)=\exp\left(-\ell/R_\lambda(y,y)\right)
\\
&&\hskip2.1cm=
\exp(-\ell\int_0^\infty \nu^{(y)}(dv)(1-{\rm e}^{-\lambda v})),
\end{eqnarray}
where $\nu^{(y)}$ is the L\'evy measure of $\tau^{(y)}$.
The assumption that the speed measure does not have atoms implies that
$\tau^{(y)}$ does not have a drift.  In case $y=0$
we write simply $L,$ $\tau$ and $\nu.$
\end{description}

\noindent
{\bf 1.3} Assuming that $X$ is started from 0 we define for $t>0$ 
\begin{equation}
\label{GD}
G_t:=\sup\{s\leq t: X_s=0\}\quad{\rm and}\quad D_t:=\inf\{s\geq t: X_s=0\}.
\end{equation}
The {\sl last exit decomposition} at a fixed time $t$ states that for $u<t<v$
\begin{eqnarray}
\label{last}
&&\nonumber
\P_0(G_t\in du, X_t\in dy, D_t\in dv)
\\
&&
\hskip3cm
=p(u;0,0)\,f_{y0}(t-u)\,f_{y0}(v-t)\,du\,dv\,m(dy).
\end{eqnarray}
In fact, this trivariate distribution is only the skeleton of a more complete body of processes:
\begin{eqnarray}
\label{last1}
\{X_u: u\leq G_t\},\quad \{X_{G_t+v}: v\leq t-G_t\},\quad {\rm and}\quad \{X_{t+v}: v\leq D_t-t\}
\end{eqnarray}
the distributions of which we now characterize following  Salminen \cite{salminen97}. 
For general approaches; see
 Getoor and Sharpe \cite{getoorsharpe73a},  \cite{getoorsharpe73b}, and 
Maisonneuve \cite{maisonneuve75}.
 
Let $x,y\in\R_+$ and $u>0$ be given. Denote by $(X^{x,u,y},\P_{x,u,y})$ the diffusion bridge from $x$ to $y$ of length $u$ 
constructed from $X,$ i.e., the measure $\P_{x,u,y}$ governing $X^{x,u,y}$  is the conditional measure associated with $X$ 
started from $x$ and conditioned to be at $y$ at time $u.$  The bridge  $X^{x,u,y}$ is a strong
non-time-homogeneous Markov process defined on the time axis $[0,u).$ 
For the first component in (\ref{last1}),  we have conditionally on $G_t=u$  
\begin{eqnarray}
\label{last2}
\{X_s: 0\leq s<G_t\}\,\rr\,\{X^{0,u,0}_s: 0\leq s<u\}
\end{eqnarray}

For the second component in (\ref{last1}) consider the process $\widehat X^{(y)}$ as introduced in (iv) above with $y=0.$ 
We write simply $\widehat X$ instead of $\widehat X^{(0)}.$ For positive $x$ and $y$ let $\widehat X^{x,u,y}$ denote the bridge 
from $x$ to $y$ of length $u$ constructed, as above, from $\widehat X.$ The measure $\widehat \P_{x,u,y}$ governing $\widehat X^{x,u,y}$ 
can be extended by taking (in the weak sense)
$$
\widehat \P_{0,u,y}:=\lim_{x\downarrow 0} \widehat \P_{x,u,y}.
$$
We let  $\widehat X^{0,u,y}$ denote the process associated with $\widehat \P_{0,u,y}.$ Then, 
conditionally on $G_t=u$  and $X_t=y,$
\begin{eqnarray}
\label{last3}
 \{X_{G_t+s}: 0\leq s< t-G_t\}\,\rr\,\{\widehat X^{0,t-u,y}_s: 0\leq s<t-u\}.
\end{eqnarray}

For the final part in (\ref{last1}), by the Markov property, we have conditionally on $X_t=y$
\begin{eqnarray}
\label{last4}
 \{X_{t+s}: s< D_t-t\}\,\rr\,\{\widehat X_s: s\geq 0\},
\end{eqnarray}
where $\widehat X_0=y.$

\section{Two descriptions of the It\^o measure}
\label{sec2}

\subsection{Description due to It\^o and McKean}% (and L\'evy ?)}
\label{sec21}
We discuss the description of the It\^o measure ${\bf n}$ where the excursions are studied by 
conditioning with respect to their lifetimes. 
Let $\widehat X$ be as in section 1.3 and $\hat p(t;x,y)$ its 
transition density with respect to the speed measure, in other words,
$$
\P_x(\widehat X_t\in dy)=\P_x(X_t\in dy; t<H_0)=\hat p(t;x,y)\,m(dy).
$$
The L\'evy measure $\nu$ of $\tau$ is absolutely continuous with respect to the Lebesgue measure, 
and the density - which we also denote by $\nu$ - is given by %can be computed as follows
\begin{eqnarray}
\label{e21}
&&\hskip-.5cm
%\nonumber
\nu(v):=\nu(dv)/dv=\lim_{x\downarrow 0}\lim_{y\downarrow 0}\frac{\hat p(v;x,y)}{S(x)S(y)}
%=-\frac{\partial^2}{\partial S(x)\partial S(x)}\hat p(v;x,x)\vert_{x=0}
%\\
%&&\hskip1.1cm
=:p^\uparrow(v;0,0).
\end{eqnarray}

In Section 1.3 we have introduced the bridge $\widehat X^{x,t,y}$ and the measure
$\widehat \P_{x,t,y}$ associated with it. The family of probability
measures $\{\widehat \P_{x,t,y}: x>0,y>0\}$ is 
weakly convergent as $y\downarrow 0$ thus defining $\widehat \P_{x,t,0}$ for
all $x>0.$ Intuitively, this is the process $\widehat X$ conditioned
to hit 0 at time $t.$ Moreover, letting 
now $x\downarrow 0$ we obtain a measure which we denote by $\widehat \P_{0,t,0}$ which governs a non-time 
homogeneous Markov process $\widehat X^{0,t,0}$ starting from 0, staying positive on the time interval $(0,t)$ and ending at 0 at time $t.$

\begin{thm}\label{thm1} {\bf a.} The law of the excursion life time
  $\zeta$ under the It\^o excursion measure  ${\bf n}$
is equal to the L\'evy measure of 
the subordinator $\{\tau_\ell\}_{\ell\geq 0}$ and is given by
%\hfill\break\hfill
\begin{equation}
\label{e211}
{\bf n}(\zeta\in dv)=\nu(dv)=p^\uparrow(v;0,0)\, dv.
\end{equation}
{\bf b.} The It\^o measure can be represented as the following integral
\begin{equation}
\label{e2115}
{\bf n}(d\varepsilon)=\int_0^\infty {\bf n}(\zeta\in dv)\,\widehat \P_{0,v,0}(d\varepsilon).
\end{equation}
Moreover, the finite dimensional distributions of the excursion are characterized for 
$0<t_1<t_2<\dots<t_n$ and  $x_i>0,\ i=1,2,\dots,n$
by
\begin{eqnarray}
\label{e2116}
&&
\hskip-1cm
\nonumber
{\bf n}(\varepsilon_{t_1}\in dx_1,\varepsilon_{t_2}\in dx_2,\dots,\varepsilon_{t_n}\in dx_n)\\
&&
\hskip1cm
=
 m(dx_1)\,f_{x_{1}0}(t_1)\,\hat p(t_2-t_1;x_1,x_2)\,m(dx_2)
\\
&&
\hskip4cm
\nonumber
\times\dots
\hat p(t_n-t_{n-1};x_{n-1},x_{n})\,m(dx_{n}).
\end{eqnarray}
In particular, the excursion entrance law is given by
$$
{\bf n}(\varepsilon_t\in dx)= m(dx)\,f_{x0}(t),
$$
and it holds
\begin{equation}
\label{e212}
{\bf n}(\zeta >t)=\int_0^\infty {\bf n}(\varepsilon_t\in dx)= \int_0^\infty m(dx)\,f_{x0}(t).
\end{equation}
\end{thm}

Combining the formulas (\ref{e211}) and (\ref{e2115}) with the last
exit decomposition (\ref{last}) leads to a curious relation between
the transition densities $p$ and $p^\uparrow.$ 
\begin{proposition}
\label{prop1}
The functions $p(t;0,0)$ and $p^\uparrow(t;0,0)$ satisfy the identity
\begin{equation}
\label{marc2}
\int_0^t du\ p(u;0,0)
\int_{t-u}^\infty dv\ p^\uparrow(v;0,0) =1.
\end{equation}
\end{proposition}
\begin{proof} From (\ref{e211}) and (\ref{e212}) we may write 
$$
\int_t^\infty dv\ p^\uparrow(v;0,0)={\bf n}(\zeta >t)=\int_0^\infty m(dx)\,f_{x0}(t).
$$
Consequently, identity (\ref{marc2}) can be rewritten as 
\begin{equation}
\label{marc3}
\int_0^t du\ p(u;0,0) \int_0^\infty m(dx)\,f_{x0}(t-u)=1,
\end{equation}
but, in view of the last exit decomposition (\ref{last}), identity
(\ref{marc3}) states  that the last exit from 0 when 
starting from 0 takes place with probability 1 before $t,$ in other words,
$$
\P_0(G_t\leq t)=1,
$$
which, of course, is trivially true.
\end{proof}
\begin{remark} For another approach to (\ref{marc2}) notice that it follows from (\ref{e1}) and (\ref{e21}) 
$$
\frac 1{R_\lambda(0,0)}
=\int_0^\infty dv\ p^\uparrow(v;0,0)\left(1-{\rm e}^{-\lambda v}\right).
$$
Hence, from  the definition of the Green kernel, 
\begin{eqnarray}
\label{marc}
&&
1=
\int_0^\infty du\ {\rm e}^{-\lambda u}\,p(u;0,0)
\int_0^\infty dv\ p^\uparrow(v;0,0)(1-{\rm e}^{-\lambda v}).
\end{eqnarray}
Consequently,
\begin{eqnarray*}
&&
%\hskip-.5cm
\frac 1\lambda=
\int_0^\infty du\ {\rm e}^{-\lambda u}\,p(u;0,0)
\int_0^\infty dv\ {\rm e}^{-\lambda v} \int_v^\infty\ ds\ p^\uparrow(s;0,0)
\\
&&
\hskip.5cm
=\int_0^\infty du\,\int_0^\infty dv\, {\rm e}^{-\lambda (u+v)}\,p(u;0,0)
\, \int_v^\infty\ ds\ p^\uparrow(s;0,0),
\end{eqnarray*}
from which (\ref{marc2}) is easily deduced.
\end{remark}

\subsection{Description due to Williams}
\label{sec22}

In the approach via the lengths of
the excursions the focus is first on the time axis. In Williams'
description (see Williams \cite{williams74}, and Rogers \cite{rogers81}, \cite{rogers89}) 
the starting point of the analysis is on the space axis 
since the basic conditioning is with respect to the maximum of an excursion. 
To formulate the result, let for $\varepsilon\in E$ 
$$ 
M(\varepsilon):=\sup\{\varepsilon_t:
0<t<\zeta(\varepsilon)\}.
$$
The key element in  Williams' description is the diffusion $X^\uparrow$
obtained by conditioning $\widehat X$ not to hit 0. We use the
notation $\P^\uparrow$ and  $\E^\uparrow$ for the measure and 
the expectation 
associated with $X^\uparrow.$ To define this process rigorously set 
for a bounded $F_t\in \cC_t,\, t>0,$
\begin{eqnarray*}
&&
\E^\uparrow_x(F_t):=\lim_{a\uparrow +\infty}
\E_x(F_t\,;t<H_a\,|\, H_a<H_0) %\wedge H_0\,|\, H_a<H_0) 
\\
&&
\hskip1.3cm
=
\lim_{a\uparrow +\infty}
\frac{\E_x(F_t\,;t<H_a<H_0) }{ \P_x(H_a<H_0) }%\wedge H_0\,,\, H_a<H_0) }{ \P_x(H_a<H_0) }
\\
&&
\hskip1.3cm
=
\lim_{a\uparrow +\infty}
\frac{\E_x\left(F_t\,S(X_t)\,;t<H_a\wedge H_0\right)}{ S(x)},
\end{eqnarray*}
where the Markov property and formula (\ref{scale}) for the scale function are applied. The monotone convergence theorem 
yields 
\begin{eqnarray*}
&&
\E^\uparrow_x(F_t)=
\frac {1}{S(x)}\,\P_x\left(F_t\,S(X_t)\,;t<H_0\right), 
\end{eqnarray*}
in other words, the desired conditioning
is realized as Doob's $h$-transform of  $\widehat X$ by taking $h$ to be the scale
function of $X.$ It is easily deduced that the transition density and the speed measure associated with $X^\uparrow$ 
are given by
$$
p^\uparrow(t;x,y):= \frac{\hat p(t;x,y)}{ S(y)S(x)},\quad\,m^\uparrow(dy):=S(y)^2\,m(dy).
$$ 
We remark that the boundary point 0 is entrance-not-exit for  $X^\uparrow$ and,
therefore,  $X^\uparrow$ can be started from 0 after which it immediately
enters $(0,\infty)$ and never hits 0.

\begin{thm}\label{thm2} {\bf a.} The law of the excursion maximum $M$ under the It\^o excursion measure {\bf n} 
is given by
$$
{\bf n}(M\geq a)=\frac 1{S(a)}.
$$
{\bf b.} The It\^o excursion measure {\bf n} can be represented via
$$
{\bf n}(d\varepsilon)=\int_0^\infty {\bf n}(M\in da)\,\Q^{(*,a)}(d\ep),
$$ 
where $\Q^{(*,a)}$ is the distribution of two independent $X^\uparrow$ processes put back to back run from 0 until they first hit level $a.$
\end{thm}
As an illustration, we give the following formula
\begin{eqnarray*}
&&
{\bf n}(1-\exp(-\int_0^\zeta ds\, V(\ep_s)))
\\
&&
\hskip2cm
=\int_0^\infty {\bf n}(M\in da)\left(1-\left(\E^\uparrow_0(\exp(-\int_0^{H_a}du\, V(\omega_u)))\right)^2\right).
\end{eqnarray*}
If $V\geq 0,$ this quantity  is the L\'evy exponent of the subordinator 
$$
\left\{\int_0^{\tau_\ell}ds\, V(X_s): \ell\geq 0\right\},
$$
that is,
\begin{eqnarray*}
&&
\hskip-1cm
\E\left\{\exp\left(-\alpha\,\int_0^{\tau_\ell}ds\,
    V(X_s)\right)\right\}
\\
&&
\hskip2cm
=
\exp\left\{-\ell\,{\bf n}\left(1-\exp\left(-\alpha \int_0^{\zeta}ds\,
    V(\varepsilon_s)\right)\right)\right\}.
\end{eqnarray*}

Comparing the descriptions of the It\^o excursion measure in Theorem
\ref{thm1} (in particular formula (\ref{e2116})) and in Theorem \ref{thm2} hints
that the processes $\widehat X$ and $X^\uparrow$ have, in addition to
conditioning relationship, also a time reversal relationship. This is
due to Williams \cite{williams74}, who particularized to the case of diffusions the general 
time reversal result, obtained by Nagasawa \cite{nagasawa64}. See also \cite{revuzyor01} p. 313, and
\cite{borodinsalminen02} p. 35.
%See also \cite{sharpe88}, \cite{salminen84}. 

\begin{proposition}
\label{timrev}
Let for a given $x>0$ 
$$
\Lambda_x:=\sup\{t\,:\, \omega(t)=x\}
$$
denote the last exit time from $x.$ Then 
\begin{eqnarray}
\label{will}
\{\widehat X_s: 0\leq s<H_0\}\,\rr\,\{X^{\uparrow}_{\Lambda_x-s}: 0\leq s<\Lambda_x\},
\end{eqnarray}
where $\widehat X_0=x$ and $X^{\uparrow}_0=0.$
\end{proposition}

\section{Stationary excursions; Bismut's description}
\label{sec3.1}
Consider the diffusion $X$ with the time parameter $t$  taking values in
the whole of $\R.$ In the case $m(\R_+)<\infty$ the measure governing
$X$ can be normalized to be a probability measure. Indeed, in this
case the
distribution of $X_t$ is for every $t\in\R$ defined to be 
$$
\P(X_t\in dx)= m(dx)/m(\R_+)=:\widehat m(dx).
$$
Recall from (\ref{GD}) the definitions of $G_t$
and $D_t,$ and introduce also $\Delta_t:=D_t-G_t.$ 

\begin{thm}\label{thm31} Assume that  $m(\R_+)<\infty.$ 
Then the joint distribution of  $t-G_t$ and $D_t-t$ is given by 
\begin{eqnarray}
&&
\hskip-1cm
\nonumber
\P(t-G_t\in du, D_t-t\in dv)/du dv=
\int_0^\infty \widehat m(dy)\,
f_{y0}(u)\,f_{y0}(v)
\\
&&\hskip4.8cm
\nonumber
=
\nu(u+v)/m(\R_+).
\end{eqnarray} 
Consequently, for $\Delta_t$ it holds
\begin{equation}
\label{my0}
\P(\Delta_t\in du)/du= u\,\nu(u)/m(\R_+).
\end{equation}
Moreover, the law of the process 
$\{X_{{G_{t}}+v}:v\leq \Delta_t\}$ is given by
\begin{equation}
\label{my1}
\zeta(\varepsilon)\, {\bf n}(d\varepsilon)/m(\R_+),
\end{equation}
where ${\bf n}(d\varepsilon)$ is the It\^o measure as introduced in Theorem 
\ref{thm1} and \ref{thm2} and $\zeta$ denotes 
the length of an excursion.
\end{thm}
\begin{proof}
The density of $(t-G_t,D_t-t)$ is derived using the time reversibility of the diffusion $X,$ i.e., 
$$
\{X_t\,:\, t\in\R\}\,\rr\,\{X_{-t}\,:\, t\in\R\},
$$
and the conditional independence given $X_t.$ The fact that the
density can be expressed via the density of the L\'evy measure is
stated (and proved) in Proposition \ref{jD} below, 
see formulas (\ref{n10}) and  (\ref{n105}). To compute the distribution of $\Delta_t$ is elementary from the 
joint distribution of $t-G_t$ and $D_t-t.$ For these results, we refer
also to Kozlova and Salminen  \cite{kozlovasalminen05}. The statement concerning the law of
$\{X_{{G_{t}}+v}:v\leq \Delta_t\}$ 
has been proved in Pitman \cite{pitman86} (see  Theorem p. 290
point (iii) and the formulation for excursions on p. 293 and 294) --
all that remains for us to do is to find the right normalization constant, but
this is fairly obvious, e.g., from the density of $\Delta_t.$
\end{proof}

If $m(\R_+)=\infty$ the measure associated with $X$ is still
well-defined but ``only'' $\sigma$-finite. In this case, the
distribution of $X_t$ is plainly taken to be $m.$
From (\ref{my1}) it is seen that we are faced with a representation 
of the It\^o measure via stationary excursions valid in both cases
$m(\R_+)<\infty$ and $m(\R_+)=\infty.$ We focus now on this
representation as displayed in  (\ref{bismut}) below, and  present a
proof of the representation using the diffusion theory (this provides, of course, also a proof of
(\ref{my1})). We remark that in Pitman \cite{pitman86} a more general case 
concerning homogeneous random sets is proved, and, hence,  it
seems worthwhile to give a ``direct''  proof in the 
diffusion case.

\begin{thm} Let  $F$ be a measurable non-negative functional defined in the excursion space $E.$ Then 
up to a normalization
\begin{equation}
\label{bismut}
 {\bf n}(F(\varepsilon))=\E\left(\frac 1{\Delta_t}\ F\left(X_{G_t+s}:0\leq s\leq \Delta_t\right) \right). 
\end{equation}
In particular, the process $\{X_{G_t+s}:0\leq s\leq \Delta_t\}$ conditionally on $\Delta_t=v$ is identical in law 
with the excursion bridge $\widehat X^{0,v,0}$ as introduced in Section 2.1.
\end{thm}

\begin{proof} Without loss of generality, we take $t=0.$ From (\ref{my0}) we have 
$$
{\bf n}(f(\zeta))=\int_0^\infty f(a)\,\nu(a)\, da =\E\left(\frac 1{\Delta_0}\ f(\Delta_0)\right).
$$
Therefore, it is enough (cf. Theorem \ref{thm1}) to prove that 
%simplying the notation by introducing $Y_s:=X_{G_0+s}$ for
%$0\leq s\leq \Delta_0,$ 
\begin{equation}
\label{bis0}
{\bf n}(F(\varepsilon)\,|\,\zeta=u)=
\E\left(F\left(X_{G_0+s}:0\leq s\leq \Delta_0\right)\,|\, \Delta_0=u \right). 
\end{equation}
Define for $0\leq s_1<s_2<\dots<s_n\leq u$%\Delta_0$ 
$$
A_{1,n}:=\{X_{G_0+s_1}\in dy_1,\dots, X_{G_0+s_n}\in dy_n\},
$$
and consider
\begin{eqnarray*}
&&
\E\left(A_{1,n}\,|\, \Delta_0=u \right)=\int_{y=0}^\infty\int_{v=0}^u 
\E\left(A_{1,n},\,-G_0\in dv,\, X_0\in dy\,|\, \Delta_0=u \right)
\\
&&
\hskip2.9cm
= 
\int_{y=0}^\infty\int_{v=0}^u 
\E\left(A_{1,n}\,|\, \Delta_0=u,\,G_0=-v,\, X_0=y \right)
\\
&&
\hskip5.9cm
\times\P\left(-G_0\in dv,\, X_0\in dy\,|\,\Delta_0=u\right).
\end{eqnarray*}
From the description of the process $X,$ the conditional
independence and the equality of the laws of the past and future given $X_0,$
and using formula (\ref{my0}) we obtain
\begin{equation}
\label{bis2}
%\hskip-.5cm
\P\left(-G_0\in dv,\, X_0\in dy\,|\,\Delta_0=u\right)=\frac
1{u\,\nu(u)}\,f_{y,0}(v)\,f_{y,0}(u-v)\, m(dy)dv.
\end{equation}
Letting $k$ be such that $-v+s_k<t<-v+s_{k+1},$ if any, we write
applying again the conditional independence
\begin{eqnarray*}
&&
\hskip-.5cm
\E\left(A_{1,n}\,|\, \Delta_0=u,\,G_0=-v,\, X_0=y \right)
\\
&&
\hskip1.5cm
=
\E\left(A_{1,k}\,A_{k+1,n}\,|\, \Delta_0=u,\,G_0=-v,\, X_0=y \right)
\\
&&
\hskip1.5cm
= 
\E\left(A_{1,k}\,|\, G_0=-v,\, X_0=y \right)
\E\left(A_{k+1,n}\,|\, D_0=u-v,\, X_0=y \right).
\end{eqnarray*}
Recall from Introduction section 1.2 (iv) the notation $\widehat X$ for the diffusion $X$ killed when it hits 0. 
As in section 1.3 of Introduction we may construct the bridge $\widehat X^{y,v,0}$ starting from $y$ having the 
length $v$ and ending at 0. We let $\widehat \P_{y,v,0}$ denote the measure associated with $\widehat X^{y,v,0}.$
With these new notations, 
\begin{eqnarray*}
&&
\hskip-.5cm
\E\left(A_{1,k}\,|\, G_0=-v,\, X_0=y \right)
\\
&&
\hskip1cm
=
\widehat \P_{y,v,0}\left(\omega_{v-s_k}\in dy_k,\dots,\omega_{v-s_1}\in dy_1\right)
\\
&&
\hskip1cm
=\frac 1{f_{y0}(v)}\,\widehat p(v-s_k;y,y_k)\,m(dy_k)\,\widehat p(s_k-s_{k-1};y_k,y_{k-1})
\,m(dy_{k-1})
\\
&&
\hskip3.5cm
\times\dots\,\widehat p(s_2-s_1;y_2,y_1)\,m(dy_1)\, f_{y_10}(s_1)
\end{eqnarray*}
and 
\begin{eqnarray*}
&&
\hskip-.5cm
\E\left(A_{k+1,n}\,|\, D_0=u-v,\, X_0=y \right)
\\
&&
\hskip1cm
=
\widehat \P_{y,v,0}\left(\omega_{s_{k+1}-v}\in dy_{k+1},\dots,\omega_{s_n-v}\in dy_n\right)
\\
&&
\hskip1cm
=\frac 1{f_{y0}(u-v)}\,\widehat p(s_{k+1}-v;y,y_{k+1})\,m(dy_{k+1})
\\
&&
\hskip2.5cm
\times
\,\widehat p(s_{k+2}-s_{k+1};y_{k+1},y_{k+2})
\,m(dy_{k+2})
\\
&&
\hskip3.8cm
\times
\,\dots\,\widehat p(s_n-s_{n-1};y_{n-1},y_n)\,m(dy_n)\, f_{y_n0}(u-s_n).
\end{eqnarray*}
Using now (\ref{bis2}) and formulas above we have after some rearranging and applying the symmetry of the transition 
density $\widehat p$ 
\begin{eqnarray*}
&&
\E\left(A_{1,n}\,|\, \Delta_0=u \right)
\\
&&
\hskip1cm
=
\frac 1{u\,\nu(u)}\,m(dy_1)\, f_{y_10}(s_1)
\widehat p(s_2-s_1;y_1,y_2)\,m(dy_2)
\\
&&
\hskip2cm
\times\dots
\,\int_{0}^udv\, \int_{0}^\infty m(dy)\,\widehat p(v-s_k;y_k,y)\,
\widehat p(s_{k+1}-v;y,y_{k+1})%\,m(dy_{k+1})
\\
&&
\hskip3cm
\times
\dots
\,\widehat p(s_n-s_{n-1};y_{n-1},y_n)\,m(dy_n)\, f_{y_n0}(u-s_n).
\end{eqnarray*}
Performimg the integration yields
\begin{eqnarray*}
&&
\E\left(A_{1,n}\,|\, \Delta_0=u \right)
\\
&&
\hskip1cm
=
\frac 1{\nu(u)}\,m(dy_1)\, f_{y_10}(s_1)
\widehat p(s_2-s_1;y_1,y_2)\,m(dy_2)
\\
&&
\hskip2cm
\times\dots
\,\widehat p(s_n-s_{n-1};y_{n-1},y_n)\,m(dy_n)\, f_{y_n0}(u-s_n),
\end{eqnarray*}
and this means that (\ref{bis0}) holds completing the proof.
\end{proof}
\begin{remark}
The formula (\ref{bismut}) was derived for Brownian motion 
by Bismut \cite{bismut85}. The connection with the Palm measure and
stationary processes is discussed in Pitman \cite{pitman86}. In fact, Bismut describes in the Brownian case 
the law of the process 
$\{X_{G_t+s}:0\leq s\leq \Delta_t\}$ in terms of two independent 3-dimensional Bessel processes 
started from 0 and killed 
at the last exit time from an independent level distributed according to the Lebesgue measure (see
 \cite{bismut85} and \cite{revuzyor01} for details).
\end{remark}
%%h

\section{On the excursion straddling an independent exponential time}
\label{sec4}

In the literature one can find several papers devoted to the properties of 
excursions straddling a fixed time $t;$ first of all, L\'evy's fundamental paper \cite{levy39}, which 
contains a lot about the zero set of Brownian motion, its (inverse) local time, excursions, and so on. See also
Chung \cite{chung76} starting from L\'evy's paper \cite{levy39}, Durret and Iglehart \cite{durrettiglehart77}, 
and Getoor and Sharpe \cite{getoorsharpe79}, \cite{getoorsharpe82}. In fact,
the last exit decomposition (\ref{last}) lies in the heart of these
studies (see  Getoor and Sharpe \cite{getoorsharpe73a}, \cite{getoorsharpe73b}). 
However, it seems to us that excursions straddling an exponential time are not so
much analyzed. Here we make some remarks on this subject. 

Let $T$ be an exponentially distributed random variable 
with parameter $\ \alpha>0,$ independent of $X,$ and define
$$
G_T:=\sup\{s\leq T: X_s=0\},\quad\quad D_T:=\inf\{s\geq T: X_s=0\},
$$
and
$$
\Delta_T:=D_T-G_T.
$$
The L\'evy exponent of the inverse local time at 0 is denoted by
$\Phi,$ in other words,
$$
\E_0\left(\exp(-\lambda \tau_\ell)\right)=\exp\left(-\ell\,\Phi(\lambda)\right)
$$
Recall the relation (cf. (\ref{e1}) with $y=0$)
\begin{equation}
\label{BER1}
\Phi(\lambda)\,R_\lambda(0,0)=1.
\end{equation}

\subsection{Last exit decomposition at $T$}
\label{ssec41}
We begin by discussing the last exit decomposition at the exponential time $T.$ %(having mean $1/\alpha.$ 

\begin{thm} \label{thmarc} {\bf (i)} The processes 
$$\{X_u:u\leq G_{T}\}\quad{\rm and}\quad   
\{X_{{G_{T}}+v}:v\leq  \Delta_{T}\}
$$ 
are independent.
\hfill\break\hfill
{\bf (ii)} The law of  $\{X_u:u\leq G_{T}\}$ may be described as follows:
\begin{description}
\item{(a)} $L_{T}:= L_{G_{T}}$ is exponentially distributed with mean $1/\Phi(\alpha),$%=R_\alpha(0,0),$
\item{(b)} The process  $\{X_u:u\leq G_{T}\}$ conditionally on $L_{T}=\ell$ 
is distributed as $\{X_u:u\leq \tau_\ell\}$ under the probability
$$
\exp\left(-\alpha\,\tau_\ell+\ell\,\Phi(\alpha)\right)\, \P_0.
$$
\end{description}
\hfill\break\hfill
{\bf (iii)} The law of the process 
$\{X_{{G_{T}}+v}:v\leq \Delta_{T}\}$ is given by
\begin{equation}
\label{my}
\frac 1{\Phi(\alpha)}\left(1-{\rm e}^{-\alpha \zeta(\varepsilon)}\right) {\bf n}(d\varepsilon).
\end{equation}
where ${\bf n}(d\varepsilon)$ is the It\^o measure associated with the excursions away from 0 for $X$ and $\zeta$ denotes 
the length of an excursion.
\end{thm}
\begin{proof}
Let $F_1$ and $F_2$ be two nonnegative functionals of continuous functions and consider
\begin{eqnarray*}
&&
\E_0\left(F_1(X_u : u\leq G_T)\,F_2(X_{G_T+v} : v\leq  \Delta_{T})\right)\\
&&
\hskip1cm
=\alpha\int_0^\infty dt\ {\rm e}^{-\alpha t} \E_0\left(F_1(X_u : u\leq G_t)\,F_2(X_{G_t+v} : v\leq  \Delta_{t})\right)\\
&&
\hskip1cm
=\alpha \E_0\left(\sum_{\ell}\int_{\tau_{\ell-}}^{\tau_\ell} dt\ {\rm e}^{-\alpha t}F_1(X_u : u\leq \tau_{\ell-})
\,F_2(X_{ \tau_{\ell-}+v} : v\leq  \tau_{\ell}- \tau_{\ell-})\right)\\
&&
\hskip1cm
= \E_0\left(\int_0^\infty d\ell\ {\rm e}^{-\alpha\,\tau_\ell}F_1(X_u : u\leq \tau_{\ell})\right)\\
&&
\hskip4cm
\times
\int {\bf n}(d\varepsilon)\left(1-{\rm e}^{-\alpha\,\zeta(\varepsilon)}\right) F_2(\varepsilon_s : s\leq \zeta(\varepsilon))\\
&&
\hskip1cm
= \Phi(\alpha)\E_0\left(\int_0^\infty d\ell\ {\rm e}^{-\alpha\,\tau_\ell}F_1(X_u : u\leq \tau_{\ell})\right)\\
&&
\hskip4cm
\times
\int {\bf n}(d\varepsilon)\left(\frac{1-{\rm e}^{-\alpha\,\zeta(\varepsilon)}}{\Phi(\alpha)}\right) F_2(\varepsilon_s : 
s\leq \zeta(\varepsilon)).
\end{eqnarray*}
where the third equality is based on the properties of the Poisson random measure associated with the excursions (see 
Revuz and Yor \cite{revuzyor01} Master Formula p. 475). 
%Recall that $n$ denotes the It\^o measure, see Section 2.  
\end{proof}
\begin{remark}
\label{RM1}
Notice that letting  $\alpha\to 0$ in (\ref{my}) and using 
\begin{equation}
\label{plim}
\lim_{\alpha\to 0}\frac \alpha{\Phi(\alpha)}=\lim_{\alpha\to 0} \alpha\,R_\alpha(0,0)= 1/m(\R_+)
\end{equation}
yield the 
probability law of the excursion straddling a fixed time in the
stationary setting, cf. (\ref{my1}) in   Theorem \ref{thm31}.
\end{remark}

As a corollary of Theorem \ref{thmarc} we have the following results 
which show that after conditioning the quantities do not depend on $\alpha,$ 
and in this context $\alpha$ is entirely "contained" in $G_T$ and $\Delta_T.$
The formulas should be compared with
(\ref{last2}), (\ref{last3}), and (\ref{last4}).
The distributions of  $G_T$ and $\Delta_T$ are given, respectively, in (\ref{GT1}) and (\ref{421}) below.

\begin{corollary} For any nonnegative functionals $F_1$ and $F_2$ of continuous functions it holds
\begin{eqnarray}
\label{n4}
&&
\hskip-1.5cm
\nonumber
\E_0\left(F_1(X_u : u\leq G_T)\,|\, G_T=g\right)=\E_0\left(F_1(X_u : u\leq
  g)\,|\, X_g=0\right)
\\
&&
\hskip3.55cm
=\E_{0,g,0}\left(F_1(\omega_u : u\leq
  g)\right)
\end{eqnarray}
and
\begin{eqnarray}
\label{n5}
&&
\nonumber
\hskip-1.5cm
\E_0\left(F_2(X_{G_T+v} : v\leq \Delta_T)\,|\,  \Delta_T=h\right)
%\\
%&&
%\hskip3cm
=\E_0^\uparrow\left(F_2(X_v : v\leq h)\,|\, X_h=0\right).
\\
&&
\hskip4.3cm
=\widehat \E_{0,h,0}(F_2(\omega_s : 
s\leq h))
\end{eqnarray}
\end{corollary}
\begin{proof}  
The statement (\ref{n4}) can be obtained from the corresponding result for fixed time as presented in (\ref{last2}).
%hh
Also (\ref{n5}) can be derived from the fixed time result but we prefer to present here a proof based on the 
Master Formula. For this consider for $\delta>0$
\begin{eqnarray*}
&&
\E_0\left(F_2(X_{G_T+v} : v\leq \Delta_T)\,{\bf 1}_{\{h\leq  \Delta_T<h+\delta\}}\right)
\\
&&
\hskip.5cm
=
\int {\bf n}(d\varepsilon)\left(\frac{1-{\rm e}^{-\alpha\,\zeta(\varepsilon)}}{\Phi(\alpha)}\right) F_2(\varepsilon_s : 
s\leq \zeta(\varepsilon))\,{\bf 1}_{\{h\leq \zeta(\varepsilon)<h+\delta\}}.
\end{eqnarray*}
Using the description (\ref{e2115}) of the It\^o excursion law we obtain
\begin{eqnarray*}
&&
\E_0\left(F_2(X_{G_T+v} : v\leq  \Delta_T)\,{\bf 1}_{\{h\leq  \Delta_T<h+\delta\}}\right)
\\
&&
\hskip.5cm
=
\int_h^{h+\delta}du\, \nu(u)\left(\frac{1-{\rm e}^{-\alpha\,u}}{\Phi(\alpha)}\right) 
\,\widehat \E_{0,u,0}(F_2(\omega_s : 
s\leq u)).
\end{eqnarray*}
Applying the explicit form of the distribution of $\Delta_T$ as given
in (\ref{421}) and letting $\delta\downarrow 0$ leads to (\ref{n5}). 
\end{proof}

\subsection{On the distribution of $(G_T,D_T)$}
\label{ssec42}

In this section the distributions of $T-G_T,$ $D_T-T$ and $\Delta_T:=D_T-G_T$ are studied in detail.

\begin{proposition}
\label{jD}
The joint distribution of $T-G_T$ and $D_T-T$ is given by
\begin{eqnarray}
&&
\nonumber
\P_0(T-G_T\in du, D_T-T\in dv)
\\
&&\hskip1cm
\label{n10}
=
du dv\,\alpha\, R_\alpha(0,0)\,{\rm e}^{-\alpha u}\,\int_0^\infty m(dy)\,
f_{y0}(u)\,f_{y0}(v).
\\
&&\hskip1cm
\label{n105}
=
\frac{\alpha\,{\rm e}^{-\alpha u}\,\nu(u+v)}{ \Phi(\alpha)}\,du dv.
\end{eqnarray} 
In particular,
\begin{equation}
\label{n11}
\nu(u+v)=\int_0^\infty m(dy)\,
f_{y0}(u)\,f_{y0}(v).
\end{equation}
\end{proposition}

\begin{proof} From (\ref{last}),
\begin{eqnarray*}
&&
\P_0(t-G_t\in du, D_t-t\in dv)
\\
&&\hskip1cm
=du dv\,p(t-u;0,0)\,{\bf 1}_{\{u\leq t, v\geq 0\}} \,\int_0^\infty m(dy)\,f_{y0}(u)\,f_{y0}(v),
\end{eqnarray*}
and, hence,
\begin{eqnarray*}
&&
\P_0(T-G_T\in du, D_T-T\in dv)
\\
&&\hskip1cm
=
du dv\,\alpha\,\int_u^\infty dt\, {\rm e}^{-\alpha t}\, p(t-u;0,0)\,\int_0^\infty m(dy)\,f_{y0}(u)\,f_{y0}(v), 
\end{eqnarray*}
from which (\ref{n10}) follows. To derive (\ref{n105}), we
apply again the Master Formula (see Revuz and Yor
\cite{revuzyor01} p. 475). For this, let $(u,v)\mapsto\varphi(u,v)$ be
a non-negative and Borel measurable function and define
$$
Q(\varphi):=\E_0(\varphi(T-G_T,D_T-T)).
$$
Letting $\tau_\ell$ denote the inverse of the local time $L$ at 0 we have
\begin{eqnarray*}
&&
Q(\varphi)=\E_0\left(\sum_{\ell\geq 0}\,\varphi(T-\tau_{\ell-},\tau_\ell-T)\,{\bf 1}_{\{\tau_{\ell-}<T<\tau_{\ell}\}}\right)
\\
&&\hskip1cm
=
\E_0\left(\int_{\R^2_+}\varphi(T-\tau_{\ell},z+\tau_{\ell}-T)\,
{\bf 1}_{\{\tau_{\ell}<T<\tau_{\ell}+z\}} \nu(z)\,dzd\ell\right),
\end{eqnarray*}
since $\{(\ell,\tau_{\ell}):\ell\geq 0\}$ is a Poisson point process with L\'evy measure $d\ell\,d\nu,$ and $T$ is 
independent of $\{\tau_{\ell}:{\ell\geq 0}\}.$ Apply next that $T$ is exponentially distributed 
to obtain
\begin{eqnarray*}
&&
Q(\varphi)=
\E_0\left(\int_{\R^2_+}\nu(z)dzd\ell\,\int^{\tau_{\ell}+z}_{\tau_{\ell}}dt\,\alpha \, {\rm e}^{-\alpha t}\,
\varphi(t-\tau_{\ell},z+\tau_{\ell}-t)\right)
\\
&&\hskip1cm
=
\alpha \E_0\left(\int_{\R^3_+}\nu(z) {\rm e}^{-\alpha (x+\tau_\ell)}\,
\varphi(x,z-x)\,{\bf 1}_{\{x\leq z\}}\,dxdzd\ell\right),
\end{eqnarray*}
where we have substituted $x=t-\tau_\ell.$ Furthermore, setting $y=z-x$ yields
\begin{eqnarray*}
&&
Q(\varphi)=
\alpha \E_0\left(\int_{\R^3_+}\nu(y+x) {\rm e}^{-\alpha (x+\tau_\ell)}\,
\varphi(x,y)\,dxdyd\ell\right)
\\
&&\hskip1cm
=
\alpha \E_0\left( \int_0^\infty {\rm e}^{-\alpha\, \tau_\ell}\, d\ell\right) \int_{\R^2_+}\varphi(x,y)\, {\rm e}^{-\alpha\, x}\, \nu(y+x)
\,dxdy,
\end{eqnarray*}
and (\ref{n105}) follows now easily from (\ref{e1}). The equality (\ref{n11}) is an immediate consequence of 
(\ref{n10}) and (\ref{n105}).
\end{proof} 

\begin{corollary} 
\label{cor00}
{\bf 1.}\ The densities for $T-G_T,$ $D_T-T,$ and $\Delta_T$ are given, respectively,  by
\begin{equation}
\label{nn1}
\P_0(T-G_T\in du)/du= 
\frac\alpha{\Phi(\alpha)}\,{\rm e}^{-\alpha u}\,\int_u^\infty \nu(z)dz,
\end{equation}
\begin{equation}
\label{nn2}
\P_0(D_T-T\in dv)/dv=
\frac\alpha{\Phi(\alpha)}\,{\rm e}^{\alpha v}\,\int_v^\infty \,{\rm e}^{-\alpha z}\nu(z)dz,
\end{equation}
and
\begin{equation}
\label{421}
\P_0(\Delta_T\in da)/da=\frac {(1-{\rm e}^{-\alpha
    a})\nu(a)}{\Phi(\alpha)}.
\end{equation}
{\bf 2.}\ The joint density of $T-G_T$ and $\Delta_T$ is 
\begin{equation}
\label{4211}
\P_0(T-G_T\in du,\, \Delta_T\in da)/du\,da=\frac\alpha {\Phi(\alpha)}\,
{\rm e}^{-\alpha u}\,\nu(a),\quad u\leq a.
\end{equation}
{\bf 3.}\  The density of $T-G_T$ conditionally on $\Delta_T=a$ is 
\begin{equation}
\label{4212}
\P_0(T-G_T\in du\,|\, \Delta_T=a)/du=\frac\alpha 
{1-{\rm e}^{-\alpha a}}
\,
{\rm e}^{-\alpha u},\quad u\leq a.
\end{equation}
\end{corollary}

\begin{proposition}
\label{jLT}
The joint Laplace transform of $G_T$ and $D_T$ is given by
\begin{eqnarray}
\label{LGD}
&&
\E_0\left(\exp\left(-\gamma_1 G_T-\gamma_2 D_T\right)\right)
=\frac{\Phi(\gamma_2+\alpha)-\Phi(\gamma_2)}{\Phi(\gamma_1+\gamma_2+\alpha)}.
\end{eqnarray}
In particular,
\begin{eqnarray}
\label{41}
&&
\E_0\left({\rm e}^{-\gamma \Delta_T}\right)
=\frac{\Phi(\gamma+\alpha)-\Phi(\gamma)}{\Phi(\alpha)},
\end{eqnarray}
and the random variables $G_T$ and $\Delta_T$ are independent.
The density of $G_T$ is given by
\begin{equation}
\label{GT1}
\P_0(G_T\in du)/du= \Phi(\alpha)\,{\rm e}^{-\alpha u}\,p(u;0,0).
\end{equation} 
\end{proposition}

\begin{proof}
The formula (\ref{GT1}) for the density of $G_T$ is obtained from (\ref{last}) by
integrating. The independence of $G_T$ and $\Delta_T$ follows immediately from (\ref{LGD}).
To derive the joint Laplace transform of $G_T$ and $D_T,$ consider
\begin{eqnarray*}
&&
\E_0\left(\exp\left(-\gamma_1 G_T-\gamma_2 D_T\right)\right)
\\
&&\hskip1cm
=
\int_u^vdt\,\alpha\,{\rm e}^{-\alpha t} \int\,{\rm e}^{-\gamma_1 u-\gamma_2 v}
\P_0(G_t\in du, D_t\in dv).
\end{eqnarray*}
Applying the last exit decomposition formula (\ref{last}) yields
\begin{eqnarray*}
&&
\E_0\left(\exp\left(-\gamma_1 G_T-\gamma_2 D_T\right)\right)
\\
&&\hskip1cm
=
\int_0^\infty du\,{\rm e}^{-\gamma_1 u}\,p(u;0,0)\,\int_u^\infty dv\,{\rm e}^{-\gamma_2 v}
\int_u^vdt\,\alpha\,{\rm e}^{-\alpha t}
\\
&&\hskip5cm
\times\int_0^\infty m(dy)\,f_{y0}(t-u)\,f_{y0}(v-t)
\\
&&\hskip1cm
=
\int_0^\infty du\,{\rm e}^{-\gamma_1 u}\,p(u;0,0)\,\int_0^\infty da\,{\rm e}^{-\gamma_2(a+u)}
\int_u^{a+u}dt\,\alpha\,{\rm e}^{-\alpha t}
p(u;0,0)
\\
&&\hskip5cm
\times\int_0^\infty m(dy)\,f_{y0}(t-u)\,f_{y0}(a+u-t)
\\
&&\hskip1cm
=
\int_0^\infty du\,{\rm e}^{-(\gamma_1+\gamma_2) u}\,p(u;0,0)\,\int_0^\infty da\,{\rm e}^{-\gamma_2a}
\int_0^a db\,\alpha\,{\rm e}^{-\alpha(b+u) }
p(u;0,0)
\\
&&\hskip5cm
\times\int_0^\infty m(dy)\,f_{y0}(b)\,f_{y0}(a-b)
\\
&&\hskip1cm
=\alpha
\int_0^\infty du\,{\rm e}^{-(\gamma_1+\gamma_2+\alpha) u}  
p(u;0,0)\,\int_0^\infty m(dy)\,
\int_0^\infty da\,{\rm e}^{-\gamma_2a}
\\
&&\hskip5cm
\times
\int_0^a db\,{\rm e}^{-\alpha b}
f_{y0}(b)\,f_{y0}(a-b)
\\
&&\hskip1cm
=
\alpha R_{\gamma_1+\gamma_2+\alpha}(0,0)
\int_0^\infty m(dy)\,
\E_y\left({\rm e}^{-(\gamma_2+\alpha)H_0}\right)
\,\E_y\left({\rm e}^{-\gamma_2H_0}\right).
\end{eqnarray*} 
To proceed, we have 
\begin{eqnarray*}
&&
\int_0^\infty m(dy)\,
\E_y\left({\rm e}^{-(\gamma_2+\alpha)H_0}\right)
\,\E_y\left({\rm e}^{-\gamma_2 H_0}\right)
\\
&&\hskip1cm
=
\frac{1}{R_{\gamma_2+\alpha}(0,0)R_{\gamma_2}(0,0)}
\int_0^\infty m(dy)\,
R_{\gamma_2+\alpha}(y,0)R_{\gamma_2}(y,0).
\end{eqnarray*}
The integral term in this expression can be evaluated:
\begin{eqnarray*}
&&
\int_0^\infty m(dy)\,
R_{\alpha+\gamma_2}(y,0)R_{\gamma_2}(y,0)
\\
&&\hskip1cm
=
%\frac{\alpha R_{\alpha}(0,0)}{R_{\alpha+\gamma_2}(0,0)R_{\gamma_2}(0,0)}
\int_0^\infty m(dy)\,
\int_0^\infty dt\, {\rm e}^{-(\alpha+\gamma_2)t} p(t;y,0)
\int_0^\infty ds\, {\rm e}^{-\gamma_2s} p(s;y,0)
\\
&&\hskip1cm
=
%\frac{\alpha R_{\alpha}(0,0)}{R_{\alpha+\gamma_2}(0,0)R_{\gamma_2}(0,0)}
\int_0^\infty dt\, {\rm e}^{-(\alpha+\gamma_2)t} 
\int_0^\infty ds\,{\rm e}^{-\gamma_2 s} p(t+s;0,0)
\\
&&\hskip1cm
=
%\frac{\alpha R_{\alpha}(0,0)}{R_{\alpha+\gamma_2}(0,0)R_{\gamma_2}(0,0)}
\int_0^\infty dt\, {\rm e}^{-(\alpha+\gamma_2)t} 
\int_t^\infty du\, {\rm e}^{-\gamma_2(u-t)} p(u;0,0)
\\
&&\hskip1cm
=
%\frac{\alpha R_{\alpha}(0,0)}{R_{\alpha+\gamma_2}(0,0)R_{\gamma_2}(0,0)}
\int_0^\infty du\, {\rm e}^{-\gamma_2 u} \frac{1-{\rm e}^{-\alpha u}}\alpha
p(u;0,0),
\\
&&\hskip1cm
=\frac 1\alpha\left(R_{\gamma_2}(0,0)-{R_{\gamma_2+\alpha}(0,0)}\right).
\end{eqnarray*}
where the Chapman-Kolmogorov equation and the symmetry of the
transition density $p$ is applied, and by (\ref{BER1}) this completes the proof.
\end{proof}

\begin{remark}
\label{RM2}
{\bf 1.}\ From Proposition \ref{jD} it is seen that the density of $\Delta_T$
can also be written in the form
$$
P_0(\Delta_T\in da)/da=
\frac\alpha {\Phi(\alpha)}%R_{\alpha}(0,0)
\int_0^\infty m(dy)\,
\int_0^a db\,{\rm e}^{-\alpha b}
f_{y0}(b)\,f_{y0}(a-b),
$$
which taking into account (\ref{421}) leads to the identity 
$$
\frac{(1-{\rm e}^{-\alpha a})}\alpha\,\nu(a)
 =\int_0^\infty m(dy)\,
\int_0^a db\,{\rm e}^{-\alpha b}
f_{y0}(b)\,f_{y0}(a-b).
$$
Let here $\alpha\to 0$ to obtain
\begin{equation}
\label{4202}
\nu(a)=
\int_0^\infty m(dy)\,
\int_0^a \frac{db}a\,
f_{y0}(b)\,f_{y0}(a-b).
\end{equation}
It is interesting to compare this expression with the following one obtained from (\ref{n11})
\begin{equation}
\label{422}
\nu(a)=
\int_0^\infty m(dy)\,
f_{y0}(b)\,f_{y0}(a-b).
\end{equation}
The fact that the right hand sides of (\ref{4202}) and (\ref{422}) 
do not depend on $b$ can also be explained via the 
Chapman-Kolmogorov equation. 
\hfill\break\hfill
{\bf 2.}\  We may study distributions associated with $G_t,$ $D_t$ and
$\Delta_t$ in the stationary case, i.e., if $m(\R_+)<\infty,$  by letting $\alpha\to 0,$ as
observed in Remark \ref{RM1}. From Proposition \ref{jD} and Corollary
\ref{cor00} we deduce the following results take $t=0$:
$$
\P(-G_0\in du, D_0\in dv)/dudv=\frac 1{m(\R_+)}\,\nu(u+v).
$$
$$
\P(-G_0\in du)/du=\P(D_0\in du)/du=\frac {1}{m(\R_+)}\,\int_u^\infty \nu(v)\,dv,
$$
$$
\P(\Delta_0\in da)/da=\frac {1}{m(\R_+)}\,a\,\nu(a).
$$
Moreover, letting $Z_T:=(T-G_T)/\Delta_T$ then
$\left(Z_T,\Delta_T\right)$ converges in distribution as $\alpha\to 0$
to $(U,\Delta),$ where $U$ and $\Delta$ are independent with $U$
uniformly distributed on $(0,1)$ and $\Delta$ is distributed as $\Delta_0$ (cf. Theorem \ref{thm31}).
\end{remark}

\subsection{Infinite divisibility}
\label{ssec43}

In the paper by Bertoin et al. \cite{bertoinetal06} it is proved that the distribution of $\Delta_T$ for a Bessel process 
with dimension $d=2(1-\alpha), 0<\alpha<1,$ is infinitely divisible
(in fact, self-decomposable) and the L\'evy measure associated 
with this distribution is computed. In this section we show that the distribution of $\Delta_T$ is infinitely divisible 
in general, i.e., for all regular and recurrent diffusions. Moreover, we also prove that the distributions of 
$T-G_T$ and $D_T-T$ have this property. The key to these results
is the Krein representation of the density of the L\'evy measure $\nu$
 (see Knight \cite{knight81}, Kent \cite{kent82}, 
K\"uchler and Salminen \cite{kuchlersalminen89}, and, in general on Krein's theory of strings,   
Kotani and Watanabe \cite{kotaniwatanabe81}, Dym and McKean \cite{dymmckean76}) according to which
\begin{equation}
\label{429}
\nu(a)=\int_0^\infty {\rm e}^{-az} M(dz),
\end{equation}
where the measure $M$ has the properties
$$
\int_0^\infty \frac{M(dz)}{z(z+1)}<\infty\quad {\rm and}\quad \int_0^\infty \frac{M(dz)}{z}=\infty.
$$.   
\begin{thm}
\label{delta}
The distributions of $T-G_T,$ $D_T-T$ and $\Delta_T$ are infinitely divisible.
\end{thm}
\begin{proof}
As seen from   (\ref{nn1}), (\ref{nn2}), and (\ref{421}), the intrinsic
term in  the densities of  $T-G_T,$ $D_T-T$ and $\Delta_T$ is the
density $\nu(a)$ of the L\'evy measure of the inverse local time at
0. We consider first the distribution of $T-G_T.$ Applying the Krein representation 
(\ref{429}) in (\ref{nn1}) yields
\begin{eqnarray*}
&&
\hskip-1cm
\P_0(T-G_T\in du)/du
=\alpha\, R_\alpha(0,0)\,{\rm e}^{-\alpha u}\, \int_u^\infty da\,\int_0^\infty M(dz)\, {\rm e}^{-az}
\\
&&
\hskip2.7cm
=\frac \alpha{\Phi(\alpha)}\,{\rm e}^{-\alpha u}\, \int_0^\infty \frac{M(dz)}z\, {\rm e}^{-uz}
\\
&&
\hskip2.7cm
=\frac \alpha{\Phi(\alpha)}\, \int_0^\infty \frac{M(dz)}z\, {\rm e}^{-(\alpha+z)u} 
\\
&&
\hskip2.7cm
=\int_0^\infty\,(\alpha+z)\, {\rm e}^{-(\alpha+z)u}\,\widehat M_\alpha(dz),
\end{eqnarray*}
with
\begin{equation}
\label{451}
\widehat
M_\alpha(dz)= \frac \alpha{\Phi(\alpha)}\, \frac { M(dz)}
{z(\alpha+z)}.
\end{equation}
The claim of the theorem follows now from the fact that 
mixtures of exponential distributions are infinitely divisible.
(see Bondesson \cite{bondesson92}). For $D_T-T$ we compute similarly from  (\ref{nn2}) via the Krein representation  
\begin{eqnarray*}
&&
\hskip-1cm
\P_0(D_T-T\in dv)/dv
=\frac{\alpha}{\Phi(\alpha)}\,{\rm e}^{\alpha v}\,\int_v^\infty \,{\rm e}^{-\alpha a}\nu(a)da.
\\
&&
\hskip2.7cm
=\int_0^\infty z\,{\rm e}^{-zv}\,\widehat M_\alpha(dz).
\end{eqnarray*}
To analyze the distribution of $\Delta_T$ we use the Krein representation in (\ref{421}) to obtain 
\begin{equation}
\label{430}
P_0(\Delta_T\in da)/da=\frac 1{\Phi(\alpha)}\,\int_0^\infty 
\left({\rm e}^{-za}-{\rm e}^{-(\alpha+z) a}\right)\, M(dz).
\end{equation}
Notice that for $a\geq 0$ 
$$
f(a;z,\alpha)=\frac{z(\alpha +z)}{\alpha}\left({\rm e}^{-za}-{\rm e}^{-(\alpha+z) a}\right)
$$
is a probability density as a function of $a.$ In fact, letting $T_1$ and $T_2$ be two independent exponentially 
distributed random variables, with respective parameters $z$ and $\alpha+z,$ then the sum 
 $T_1+T_2$ has the density $f(a;z,\alpha).$ In particular, the distribution of
 $T_1+T_2$ is a gamma convolution (which, by definition, is the law of finite sum of independent gamma variables). 
Next we notice that 
letting
$$
\Pi_{z,\alpha}(dx):=\frac{z(\alpha+z)}{\alpha}\,x^{-2}\,dx,\quad z<x<\alpha+z
$$
we may represent the distribution of $T_1+T_2$ as a mixture of Gamma(2)-distributions as follows
\begin{equation}
\label{431}
f(a;z,\alpha)=\int_0^\infty x^2\,a\,{\rm e}^{-xa}\,\Pi_{z,\alpha}(dx).
\end{equation}
Combining the representation (\ref{431}) with (\ref{430}) yields 
\begin{eqnarray}
\label{44}
&&
\nonumber
P_0(\Delta_T\in da)/da=\frac \alpha{\Phi(\alpha)}\,\int_0^\infty 
\frac{f(a;z,\alpha)}{z(\alpha+z)} \, M(dz)
\\
&&\hskip2.9cm
= \int_0^\infty x^2\,a\,{\rm e}^{-xa}\,\widehat\Pi_{\alpha}(dx),
\end{eqnarray}
where $\widehat\Pi_{\alpha}$ is a probability measure on $\R_+$ given
for any Borel set $A$ in $\R_+$ by
\begin{equation}
\label{45}
\widehat\Pi_{\alpha}(A)=\int_0^\infty \widehat
M_\alpha(dz)\Pi_{z,\alpha}(A).
\end{equation}
The claim that the distribution of $\Delta_T$ is infinitely divisible follows now from (\ref{44}) by
evoking the result that mixtures of Gamma(2)-distributions are
infinitely divisible (see Kristiansen \cite{kristiansen94}).
\end{proof}
\begin{remark} 
%t 
{\bf 1.} Recall from Bondesson \cite{bondesson92} that a probability distribution $F$ on
$\R_+$ is called a generalized gamma convolution (GGC) if its Laplace
transform can be written as
\begin{equation}
\label{43}
\int_0^\infty {\rm e}^{-sa}F(da)=\exp\left(-\mu s+ \int_0^\infty\log\left(\frac t{t+s}\right) U(dt)\right),
\end{equation}
where $\mu\geq 0$ and $U$ is a measure on  $(0,\infty)$ satisfying
$$
\int_{(0,1]} |\log t| U(dt)<\infty\quad {\rm and}\quad \int_{(1,\infty)} \frac{U(dt)}{t}<\infty.
$$    
It is known see \cite{bondesson92} Theorem 4.1.1 p. 49 that if $\beta$ is the total mass of $U$ 
then the distribution $F$ in (\ref{43}) is a mixture of Gamma($\beta$)-distributions  
\hfill \break\hfill\noindent
{\bf 2.} The distribution of the length $\Delta_t$ of an excursion
  straddling a fixed time $t$ for a stationary diffusion (with
  stationary probability distribution) is given in Theorem \ref{thm31} (\ref{my0})
as
$$ 
\P(\Delta_t\in da)=\frac {a\,\nu(a)}{m(\R_+)}\,da.
$$ 
 Also in this case the distribution of $\Delta_t$ is a
mixture of Gamma(2)-distributions and, hence, it is infinitely divisible. In
fact, 
$$
 \P(\Delta_t\in da)/da=\int_0^\infty z^2\,a\, {\rm e}^{-za}\widetilde
 M(dz).
$$
where the probability measure $\widetilde M$ is given in
terms of the Krein measure $M$ via
$$
\widetilde M(dz)=M(dz)/(m(\R_{+})z^2).
$$
\end{remark}

%h

\section{Case study: Ornstein-Uhlenbeck processes}
\label{sec5}

In this section we give some explicit formulas for excursions from 0 to 0 associated with Ornstein-Uhlenbeck
processes. It is possible to obtain 
such formulas due to the symmetry of  the Ornstein-Uhlenbeck
process around 0. Analogous results for excursions from an arbitrary point $x$ to $x$ are 
less tractable. 

\subsection{Basics}
\label{sec51}

Let $U$ denote the Ornstein-Uhlenbeck diffusion with parameter $\gamma>0,$ i.e., $U$ is the 
solution of the SDE
$$
dU_t=dB_t-\gamma U_tdt\quad {\rm with}\quad U_0=u,
$$
and most of the time, but not always, we take $u=0.$ Recall that the speed measure and the scale function of $U$ can be taken to be 
$$
m(dx):= 2 \, {\rm e}^{-\gamma x^2}\, dx\quad {\rm and}\quad S(x):=\int_0^x{\rm e}^{\gamma y^2}\, dy,
$$
respectively. Moreover, see \cite{borodinsalminen02} p. 137, the Green kernel of Ornstein-Uhlenbeck
process with respect to the speed measure is given for $x\geq y$ by
\begin{eqnarray*}
&&
R_\lambda(x,y)=\frac{\Gamma(\lambda/\gamma)}{2\sqrt{\gamma\pi}}\, \exp\left(\frac{\gamma x^2}2\right)\,D_
{-\lambda/\gamma}\left(x\sqrt{2\gamma}\right)
\\   
&&\hskip4cm
\times \exp\left(\frac{\gamma y^2}2\right)\,D_
 {-\lambda/\gamma}\left(-y\sqrt{2\gamma}\right),
\end{eqnarray*}
where $D$ denotes the parabolic cylinder function.
In particular, since 
$$
D_{-\lambda/\gamma}(0)=
\sqrt{\pi}\left(2^{\lambda/(2\gamma)}\,\Gamma\left(( \lambda+\gamma)/(2\gamma)\right)\right)^{-1},
$$
we have, after some manipulations, 
$$
R_\lambda(0,0)=
\frac{\sqrt{\pi}\,\Gamma(\lambda/\gamma)}
{2\,\sqrt{\gamma}} 
\left(2^{\lambda/(2\gamma)}\,
\Gamma\left(( \lambda+\gamma)/(2\gamma)\right)\right)^{-2}.
$$
Consequently, using the formula
$$
\Gamma(x)=\frac{2^{\,x-1}}{\sqrt\pi}\,\Gamma((x +1)/2)\,\Gamma(x/2)
$$
we obtain
\begin{equation}
\label{ber}
R_\lambda(0,0)=\frac 1{\Phi(\lambda)}=\frac{\Gamma(\lambda/(2\gamma))}
 {4\,\Gamma((\lambda +\gamma)/(2\gamma))}.
\end{equation}

We remind also that 
$U$ can be represented as the deterministic time change (Doob's transformation) of Brownian motion via
$$  
U_t={\rm e}^{-\gamma t}\left(u+\beta_{a_t}\right),
$$
where $\beta$ is a standard Brownian motion and $a_t:=({\rm e}^{2\gamma t}-1)/2\gamma$ 
(see Doob \cite{doob42}).
% , see also Breiman \cite{breiman68}). 

\subsection{Killed Ornstein-Uhlenbeck processes}
\label{sec52}

We consider now the Ornstein-Uhlenbeck process killed at the first hitting time of 0, and denote this 
process by $\widehat U.$ Let $Y$ be the diffusion on $\R_+$ satisfying the SDE
$$
dY_t=dB_t+\left(\frac 1{Y_t}-\gamma Y_t\right)dt,\quad Y_0=y>0.
$$
Recall that $Y$ may be described as the radial part of the three-dimensional 
Ornstein-Uhlenbeck process. In \cite{borodinsalminen02} p. 138 the basic properties of such processes are presented. 
In particular, we record that 0 is an entrance-not-exit boundary point and the process is positively recurrent 
its stationary distribution being the Maxwell distribution, i.e.,
the distribution with the density proportional to the speed measure of $Y,$ that is,
$$
m^Y(dx):= 2 x^2\, {\rm e}^{-\gamma x^2}\, dx,\quad x>0.
$$
We remark that there is a misprint in \cite{borodinsalminen02} 
p. 139; the stationary distribution in the general case is not a $\chi^2$-distribution 
but a generalization of the Maxwell distribution.  
The transition density of $Y$ with respect to its speed measure $m^Y$ is
$$
p^Y(t;x,y)=\frac{\sqrt{\gamma}\,{\rm e}^{3\gamma t/2}}{\sqrt{2\pi\,\sh(\gamma t)}\,xy}
\,\exp\left( -\frac{\gamma {\rm e}^{-\gamma t}(x^2+y^2)}{2\,\sh(\gamma t)}\right)\,
\sh\left(\frac{\gamma xy}{\sh(\gamma t)}\right)
$$
and can be computed from the transition density of a Bessel process using Doob's transform (for an approach 
via inverting the Laplace transform see Giorno et al. \cite{giornoetal86}).
In Salminen \cite{salminen84b} it is proved that 
\begin{eqnarray}
\label{den}
&&\nonumber
\P_x(\widehat U_t\in dy)=\P_x(U_t\in dy, t<H_0)
\\
&&
\hskip2.2cm
={\rm e}^{-\gamma t}\,p^Y(t;x,y)\,\frac{\varphi_\gamma(y)}{\varphi_\gamma(x)}\, m^Y(dy),
\end{eqnarray}
where $\varphi_\gamma(x)=1/x$ is the unique (up to multiplicative constants) decreasing positive solution of 
the ODE associated with $Y$ killed at rate $\gamma$:
$$
\frac 12 u''(x)+\left(\frac 1x-\gamma x\right)u'(x)=\gamma u(x).
$$
From (\ref{den}) we obtain 

\begin{proposition}
\label{killedOU}
The transition density (with respect to its speed measure $m$) of the Ornstein-Uhlenbeck killed at the first hitting time of 0
is given by
\begin{equation}
\label{hatp}
\hat p(t;x,y)
=\frac{\sqrt{\gamma}\,{\rm e}^{\gamma t/2}}{\sqrt{2\pi\sh(\gamma t)}}
\,\exp\left( -\frac{\gamma {\rm e}^{-\gamma t}(x^2+y^2)}{2\sh(\gamma t)}\right)\,
\sh\left(\frac{\gamma xy}{\sh(\gamma t)}\right).
\end{equation}
\end{proposition}
%\noindent 

Combining the expression of the transition density in (\ref{hatp}) with formula (\ref{f00}) 
yields the distribution of $H_0$ (see also Sato \cite{sato77} and Going-Jaeschke and Yor \cite{going-jaeschkeyor03}). 

\begin{proposition}
\label{hitOU}
The density of the first hitting time of 0 for the Ornstein-Uhlenbeck process $\{U_t\}$ is given by 
\begin{eqnarray}
\label{hit}
&&
\hskip-.6cm
f_{x0}(t)=
%\P_x(H_0\in dt)/dt=
\frac{\gamma^{3/2}\,x\,{\rm e}^{\gamma t/2}}{\sqrt{2\pi}(\sh(\gamma t))^{3/2}}
\,\exp\left( -\frac{\gamma {\rm e}^{-\gamma t}x^2}{2\sh(\gamma t)}\right).
\end{eqnarray} 
\end{proposition}
%\noindent

\subsection{L\'evy measure of inverse local time and densities of $\Delta_T,$ $T-G_T$ and $D_T-T$}
\label{sec53}

The density of the L\'evy measure of the inverse local time at 0 is obtained by applying 
formula (\ref{e21}) (see also Hawkes and Truman \cite{hawkestruman91}). Moreover, using (\ref{ber}) in formula (\ref{e1}) 
leads to an explicit expression for the Bernstein 
function associated with the inverse local time at
0. 
%hhh
\begin{proposition}
\label{levyOU}
The density of the L\'evy measure of the inverse local time at 0 is 
\begin{eqnarray}
\label{levden}
&&
%\hskip-.5cm
\nu(t)
=
\frac{\gamma^{3/2}\,{\rm e}^{\gamma t/2}}{\sqrt{2\pi}(\sh(\gamma t))^{3/2}}
= 
\frac{(2\gamma)^{3/2}\,{\rm e}^{2\gamma t}}{\sqrt{2\pi}\,({\rm e}^{2\gamma t}-1)^{3/2}}.
\end{eqnarray}
Let  $\{\tau_\ell : \ell\geq 0\}$ be the inverse local time at
0. Then 
$$
\E_0\left(\exp(-\lambda \tau_\ell)\right)
%=\exp\left(-\frac{\ell}{R_\lambda(0,0)}\right)
=\exp\left(-\ell\, \frac{4\,\Gamma\left((\lambda +\gamma)/2\gamma\right)}{\Gamma(\lambda/2\gamma)}\right).
$$
\end{proposition}

Next we display the distributions of  $\Delta_T,$ $T-G_T,$ and $D_T-T.$ Recall that these distributions are 
infinitely divisible and the densities are expressable via the density of the L\'evy measure, 
as stated in Corollary \ref{cor00} formulae (\ref{nn1}) and  (\ref{nn2}), and
in Theorem \ref{delta}. To simply the notation, we take $\gamma=1.$
\begin{proposition} 
\label{rest}
With $\Phi(\alpha)$ as in (\ref{ber}), 
the distributions of  $\Delta_T,$  $T-G_T$ and $D_T-T$ are given, respectively, by
\begin{equation}
\label{51}
P_0(\Delta_T\in da)/da=
\frac{1-{\rm e}^{-\alpha    a}}{\Phi(\alpha)}\,\frac 2{\sqrt{\pi}}\,{\rm e}^{2a}\,\left({\rm
    e}^{2a}-1\right)^{-3/2},
\end{equation}
\begin{equation}
\label{510}
P_0(T-G_T\in da)/da= \frac{\alpha\,{\rm e}^{-\alpha     a}}{\Phi(\alpha)}\,
%\Gamma(\alpha/2)}{4\,\Gamma((\alpha +1)/2)}\,
\frac 2{\sqrt{\pi}}\,\left({\rm
    e}^{2a}-1\right)^{-1/2},
\end{equation}
and
\begin{equation}
\label{511}
P_0(D_T-T\in da)/da= \frac{\alpha\,{\rm e}^{\alpha     a}}{\Phi(\alpha)}\,
%\Gamma(\alpha/2)}{4\,\Gamma((\alpha +1)/2)}\,
\int_a^\infty du \,{\rm e}^{-\alpha u}
\frac 2{\sqrt{\pi}}\,{\rm e}^{2u}\,\left({\rm
    e}^{2u}-1\right)^{-3/2}.
\end{equation}

\end{proposition}

\subsection{The Krein measure}
\label{sec54}

As seen in Section 4.3,  the Krein representation plays a central r\^ole in the proof of
infinite divisibility of the distributions of $T-G_T,$ $D_T-T,$ and $\Delta_T$. 
Therefore, it seems motivated to compute the
measure $M$ (cf. (\ref{429})) in this representation for Ornstein-Uhlenbeck processes.

To start with, we give the
spectral representation of the transition density of $\hat p$ of the 
Ornstein-Uhlenbeck process killed at the first hitting time of 0. 
Instead of computing from scratch, we exploit the spectral
representation for $p^Y$ (with $\gamma=1$) as presented 
in Karlin and Taylor \cite{karlintaylor81} p. 333:
\begin{equation}
\label{53}
p^Y(t;x,y)=\sum_{n=0}^\infty w^{-1}_{n,1/2}\ {{\rm e}^{-2nt}}\,L_n^{(1/2)}(x^2)\,L_n^{(1/2)}(y^2),
\end{equation}
where $\{L_n^{(1/2)}:n=0,1,2,\dots\}$ is the famly of Laguerre polynomials with parameter $1/2$
normalized via 
\begin{equation}
\label{535}
\int_0^\infty\left(L_n^{(1/2)}(x^2)\right)^2\,m^Y(dx)
= \frac{\sqrt{\pi}}2\, {{n+\frac 12\choose n}}=:w_{n,1/2}.
\end{equation}
Notice that we consider the symmetric density with respect to the
speed measure $m^Y.$ From (\ref{den}) and (\ref{53}) the spectral representation of
$\hat p$ is now obtained immediately and is given by 
\begin{equation}
\label{54}
\hat p(t;x,y)=\sum_{n=0}^\infty w^{-1}_{n,1/2}\ {{\rm e}^{-(2n+1)t}}\,xL_n^{(1/2)}(x^2)\,yL_n^{(1/2)}(y^2).
\end{equation}
The normalization (\ref{535}) coincides with the normalization in
Erdelyi et al. \cite{erdelyi53} (see formula (2) p. 188 where the notation for the
norm is $h_n$). Therefore, from  \cite{erdelyi53} formula (13)
p. 189 we have
\begin{equation}
\label{55}
L^{(1/2)}_n(0)= {{n+\frac 12\choose n}} 
\end{equation}
 and, consequently (cf. (\ref{hit})), we obtain the spectral
representation for the density of the first hitting time of 0
\begin{eqnarray}
\label{56}
&&\nonumber
f_{x0}(t)=\sum_{n=0}^\infty w^{-1}_{n,1/2}\ {{\rm e}^{-(2n+1)t}}\,xL_n^{(1/2)}(x^2)\,L_n^{(1/2)}(0).
\\
&&
\hskip1.1cm 
=
\frac 2{\sqrt{\pi}}\sum_{n=0}^\infty \, {{\rm e}^{-(2n+1)t}}\,xL_n^{(1/2)}(x^2).
\end{eqnarray}
To find the spectral representation for the density of the L\'evy
measure we apply formula (\ref{levden}) which yields
\begin{equation}
\label{57}
\nu(t)= 
\frac 2{\sqrt{\pi}}\sum_{n=0}^\infty \, {{n+\frac 12\choose n}}\,{{\rm e}^{-(2n+1)t}}.
\end{equation}
In view of (\ref{429}), we have
\begin{proposition}
\label{krein1}
The measure $M$ in the Krein representation of $\nu$ for the Ornstein-Uhlenbeck process is 
given by
$$
M(dz)=
\frac 2{\sqrt{\pi}}\sum_{n=0}^\infty \,
{{n+\frac 12\choose n}}\,  \delta_{\{2n+1\}}(dz),
$$
where $\delta$ is the Dirac measure. 
\end{proposition}

Notice that 
$$
\nu(t)=\frac 2{\sqrt{\pi}}\,{\rm e}^{-t}\,\left(1-{\rm
    e}^{-2t}\right)^{-3/2},
$$
and, hence, (\ref{57}) is obtained also from the MacLaurin expansion
of $x\mapsto (1-x)^{-3/2}$ evaluated at $x={\rm e}^{-2t}.$
\vskip1cm
\noindent
{\bf Acknowledgement.} We thank Lennart Bondesson for co-operation concerning gamma convolutions.

%References:  Pitman and Yor \cite{pitmanyor03}, 

\bigskip
\bibliographystyle{plain}
%\bibliography{andrei}
\bibliography{yor1}
\end{document}